\documentclass[review,3p,square]{elsarticle}

\biboptions{square}
\usepackage{amsmath}
\usepackage{amssymb}
\usepackage{amsthm}
\usepackage[american]{babel}
\usepackage{microtype}

\usepackage{stmaryrd}
\usepackage{bbm}
\usepackage{enumitem}
\usepackage{fancyhdr}
\usepackage{hyperref}
\hypersetup{%
  pdftitle={BSDEs},%
  pdfsubject={BSDEs},%
  pdfauthor={ShengJun FAN, Ying Hu},%
  pdfkeywords={BSDEs},%
  pdfstartview=FitH,%
  CJKbookmarks=true,%
  bookmarksnumbered=true,%
  bookmarksopen=true,%
  colorlinks=true, linkcolor=blue, urlcolor=blue, citecolor=blue, %
}

\usepackage{cleveref}
\crefname{thm}{Theorem}{Theorems}
\crefname{pro}{Proposition}{Propositions}
\crefname{lem}{Lemma}{Lemmas}
\crefname{rmk}{Remark}{Remarks}
\crefname{cor}{Corollary}{Corollaries}
\crefname{dfn}{Definition}{Definitions}
\crefname{ex}{Example}{Examples}
\crefname{section}{Section}{Sections}
\crefname{subsection}{Subsection}{Subsections}

\newcommand{\To}{\rightarrow}
\newcommand{\as}{{\rm d}\mathbb{P}\times{\rm d} t-a.e.}

\newcommand{\ps}{\mathbb{P}-a.s.}

\newcommand{\tim}{\times}

\newcommand{\F}{\mathcal{F}}

\newcommand{\N}{\mathbb{N}}

\newcommand{\wid}{\widetilde}

\newcommand{\hcal}{\mathcal{H}}

\newcommand{\R}{{\mathbb R}}

\newtheorem{thm}{Theorem}[section]
\newtheorem{lem}[thm]{Lemma}
\newtheorem{pro}[thm]{Proposition}
\newtheorem{rmk}[thm]{Remark}

\journal{ArXiv}
\begin{document}

\begin{frontmatter}

\title{{Invariant representation for generators of general time interval quadratic {BSDE}s under stochastic growth conditions}\tnoteref{found}}
\tnotetext[found]{Supported by the Fundamental Research Funds for the central universities (No.2017XKZD11).}

\author[author1]{Guangshuo Zhou}
\ead{zhou1953363@163.com}
\author[author1,author2]{Fengjiao Du}
\ead{fjdu@vip.163.com}
\author[author1]{Shengjun Fan\corref{cor}}
\ead{f\_s\_j@126.com}
\cortext[cor]{Corresponding author}
\address[author1] {School of Mathematics,
	China University of Mining and Technology,
	Jiangsu 221116, P.R. China}
\address[author2]{School of Mathematics and Statistics, Xuzhou University of Technology, Xuzhou, Jiangsu 221018, P.R. China} \vspace{-0.5cm}

\begin{abstract}
This paper is devoted to proving a general invariant representation theorem for generators of general time interval backward stochastic differential equations, where the generator $g$ has a quadratic growth in the unknown variable $z$ and satisfies some stochastic growth conditions in the unknown variable $y$. This unifies and strengthens some known results. And, a natural and innovative idea is used to prove the representation theorem.
\end{abstract}

\begin{keyword}
Backward stochastic differential equation \sep Quadratic growth \sep \\ \hspace*{1.82cm} Invariant representation \sep General time interval \sep Stochastic growth\vspace{0.2cm}

\MSC[2021] 60H10
\end{keyword}

\end{frontmatter}
\vspace{-0.4cm}

\section{Introduction}
\label{sec:1-Introduction}
\setcounter{equation}{0}

Fix a positive integer $d>0$ and an extended real number $0<T\leqslant \infty$. Denote by $\langle\cdot,\cdot\rangle$ the usual inner product in $\R^+$. Let $(\Omega, \F, \mathbb{P})$ be a complete probability space carrying a $d$-dimensional standard Brownian motion $(B_t)_{t \geqslant 0}$, starting from $B_0=0$. Let $(\F_t)_{t \geqslant 0}$ denote the natural filtration generated by $(B_t)_{t\geqslant 0}$. We assume that $\F_T=\F$ and $(\F_t)_{t\geqslant 0}$ is right-continuous and complete.

Consider the following one-dimensional backward stochastic differential equation (BSDE for short):
\begin{equation}\label{eq:1.1}
  Y_t=\xi+\int_t^T g(s,Y_s,Z_s){\rm d}s-\int_t^T \langle Z_s , {\rm d}B_s\rangle, \ \ t\in[0,T], \ \
\end{equation}
where the terminal value $\xi $ is an $\F_T$-measurable random variable, the generator $g(\omega,t,y,z):\Omega\times[0,T]\times\R\times\R^d\To\R$ is an $(\F_t)$-progressively measurable random function for each $(y,z)$. The triple $(g,T,\xi)$ is called the parameters of BSDE \eqref{eq:1.1}, which is also denoted by BSDE $(g, T, \xi$). And, a pair of $(\F_t)$-progressively measurable processes $(Y_t,Z_t)_{t\in[0,T]}$ satisfying \eqref{eq:1.1} is called a solution of BSDE $(g,T,\xi)$.

Since Pardoux and Peng \cite{PardouxPeng1990SCL} established the first existence and uniqueness result for solutions of nonlinear BSDEs under the Lipschitz assumption on the generator $g$, the BSDE theory has been well developed. As is known to all, the representation theorem for generators of BSDEs is an crucial tool to interpret the relation between generators and solutions of BSDEs, and it plays an important role in studying the properties of generators by virtue of solutions of corresponding BSDEs. Briand et al. \cite{BriandCHMP2000ECP} established the first representation theorem for generator of finite time interval BSDE \eqref{eq:1.1}: for each $(t,y,z)\in[0,T)\times \R\times\R^d$,
\begin{equation}\label{eq:1.2}
g(t,y,z)=L^2-\lim\limits_{\varepsilon\To 0^+}\frac{1}{\varepsilon}\{Y_t(g,t+\varepsilon,y+\langle z,B_{t+\varepsilon}-B_{t}\rangle)-y\},
\end{equation}
where the generator $g$ satisfies the uniformly Lipschitz condition in unknown variables $(y,z)$ and two additional assumptions. Since then, there has been a tremendous interest in studying the corresponding representation theorem for generators of BSDEs by weakening the uniformly Lipschitz condition on $g$  and/or relaxing the finite terminal time into the general case.

In particular, we would like to mention the following works among others. First of all, by virtue of Lebesgue's Lemma Jiang \cite{Jiang2005SPA}-\cite{Jiang2008AAP} successfully eliminated two additional assumptions mentioned above and extended the representation theorem obtained in \cite{BriandCHMP2000ECP}. And, by the convolution approximation technique used in Lepeltier and San Martin \cite{LepeltierSanMartin1997SPL}, Jia \cite{Jia2010SPA}, Jia and Zhang \cite{JiaZhang2013AMS} and Fan and Jiang \cite{FanJiang2010JCAM} weakened the uniformly Lipschitz condition to the uniform continuity condition and then the continuity and linear-growth condition of the generator $g$ in $(y,z)$. Furthermore, based on the existence and uniqueness results established respectively in Kobylanski \cite{Kobylanski2000AP}, Briand and Hu \cite{BriandHu2006PTRF} and Fan and Jiang \cite{FanJiang2013AMS} for BSDEs,  Ma and Yao \cite{MaYao2010SAA}, Fan et al. \cite{FanJiangXu2011EJP}, Zheng and Li \cite{ZhengLi2015SPL}-\cite{ZhengLi2018AMASES} and Xiao and Fan \cite{XiaoFan2017SPL} investigated the representation theorem of generators of BSDEs when the generator $g$ has a general growth in $y$ and/or a quadratic growth in $z$, in which the authors developed the convolution approximation technique and utilized systematically the stopping time technique and the truncation technique. Finally, we mention that all results mentioned above dealt with the finite time interval BSDEs. To the best of our knowledge, Chen and Wang \cite{ChenWang2000JAMS} first established the existence and uniqueness of solutions for infinite time interval BSDEs, and the corresponding representation theorem of generators were obtained in Zhang and Fan \cite{ZhangFan2013SPL}.

Recently, Luo and Fan \cite{LuoFan2017SD} studied the existence of bounded solutions for general time interval quadratic BSDEs under some general stochastic growth conditions. Then, a question arises naturally: under these general assumptions, does the corresponding representation theorem for generators of BSDEs still hold? This paper gives an affirmative answer of this question under some additional mild assumptions, which strengthens all the representation results mentioned above. A natural and intrinsic idea (see \cref{lem:2.5} in Section2) instead of the convolution approximation technique is used in the proof. This simplifies greatly the procedure of proof.

The remainder of the paper is organized as follows. In next section, we introduce some notations, assumptions and lemmas used in the whole paper, and in Section 3, we state the general representation theorem and prove it.

\section{Preliminaries}
\label{sec:2-Preliminaries}
\setcounter{equation}{0}

Let $|z|$ denote the Euclidean norm of $z \in \R^d$, and let $a\wedge b := {\rm min}\{a,b\}$, $a\vee b := {\rm max}\{a,b\}$, $ a^+ := a\vee 0$ and $\R_+ := [0,+\infty)$. For each predictable subset $A$ of $\Omega\times[0,T]$, let $\mathbbm{1}_A$ equal to 1 when $(\omega,t)\in A$, and 0 otherwise. By $L^\infty(\F_T)$ we denote the set of $\F_T$-measurable bounded random variables endowed with the supremum norm $\|X\|_\infty:={\rm sup}\{x:\mathbb{P}(|X|>x)>0\}$. Moreover, let $L^1([0,T];\R_+)$ represent a set of deterministic functions $h(\cdot):[0,T]\to\R_+$ satisfying
$$\int_0^T h(t){\rm d}t<+\infty$$
and $L^\infty(\Omega;L^1([0,T];\R_+))$ the set of $(\F_t)$-progressively measurable processes $f_t(\omega):\Omega\times[0,T]\to \R_+$ such that
$$\bigg\|\int_0^T f_t(\omega) {\rm d}t\bigg\|_\infty<+\infty.$$ Finally, by $\hcal^\infty(\R_+)$ we denote the set of real-valued, continuous and $(\F_t)$-progressively measurable processes $f_t(\omega):\Omega\times[0,T]\to \R_+$ such that $$\|f_\cdot\|_{\hcal^\infty}:= \mathop{{\rm ess \ sup}}_{(\omega,t)\in\Omega\times[0,T]}|f_t(\omega)|<+\infty.$$
From now on, if there is not a special illustration, all equalities and inequalities between random variable are understood to hold $\ps$
\vspace{0.1cm}

We will use the following  assumptions on the generator $g$.
\begin{enumerate}
\renewcommand{\theenumi}{(H\arabic{enumi})}
\renewcommand{\labelenumi}{\theenumi}
\item\label{A:H1} There exist a constant $\gamma > 0$ and two stochastic processes $u_., f_. \in L^{\infty}(\Omega; L^1([0,T]; \R_+))$ such that $\as$, for each $(y,z) \in \R\tim\R^d$,
$${\rm sgn}(y)g(\omega, t, y ,z) \leqslant f_t(\omega)+u_t(\omega)|y|+\gamma |z|^2.$$
\item\label{A:H2} There exists a continuous nondecreasing function $h(\cdot):\R_+ \To \R_+$ and a $\varphi_{.}\in \mathcal{S}$ such that $\as$, for each $(y,z)\in \R\tim\R^d$,
    $$|g(\omega,t,y,z)|\leqslant \varphi_t(\omega, |y|)+h(|y|)|z|^2,$$
where $\mathcal{S}$ denotes a set of functions $\varphi_t(\omega,x):\Omega \tim [0,T] \tim \R_+ \To \R_+$ satisfying the following two conditions:

$\bullet \ \as$, the function $x \To \varphi_t(\omega,x)$  is nondecreasing;

$\bullet $ for each $x \in \R_+$, $\varphi_t(\omega,x)$ is an $(\F_t)-$progressively measurable process with  $$ \int_0^T \varphi_t(\omega,x) {\rm d}t < +\infty.$$
\item\label{A:H3} $g$ is continuous in $(y,z)$, i.e., $\as$, $(y,z)\rightarrow g(\omega, t, y, z)$ is continuous.
\end{enumerate}

The following proposition is the starting point of this paper, which comes from  Lemma 2 in \cite{FanLuo2017BKMS} and Theorem 5.1 in \cite{LuoFan2017SD}. It presents existence for the bounded solutions of general time interval quadratic BSDEs and an a priori estimate of the bounded solutions.\vspace{0.1cm}

\begin{pro}\label{pro:2.1}
Suppose that $0<T\leqslant+\infty$, $\xi\in L^{\infty}(\F_T)$ and the generator g satisfies {\rm \ref{A:H1}}-{\rm \ref{A:H3}}. Then BSDE $(g,T,\xi)$ admits a solution $(Y_t,Z_t)_{t\in[0,T]}$ such that $Y_{\cdot}\in\hcal^\infty(\R_+)$ and $\mathbf{E}\big[\int_0^T |Z_t|^2 {\rm d}t\big] <+\infty$. Moreover, for each $t\in[0,T]$, we have
\begin{equation*}\label{eq:2.1}
|Y_t|\leqslant\bigg(\|\xi\|_{\infty}+\bigg\| \int_0^T f_r {\rm d}r \bigg\|_{\infty}\bigg) {\rm exp}\bigg(\bigg\|\int_0^T u_r {\rm d}r \bigg\|_{\infty}\bigg).
\end{equation*}
\end{pro}

The following Lemmas 2.2 and 2.3 can be obtained by virtue of Lebesgue's lemma, see Lemma 3.2 in \cite{ZhangFan2013SPL} and Lemma 3.3 in \cite{ZhengLi2018AMASES} for more details.\vspace{0.1cm}

\begin{lem}\label{lem:2.2}
Let $0<T\leqslant+\infty$ and $(f_t)_{t\in[0,T\wedge m]}\in L^1([0, T\wedge m]; \R_+)$ for all $m\in\N$. Then we have
$$
  f_t=\lim\limits_{\varepsilon\To 0^+}\frac{1}{\varepsilon}\int_t^{t+\varepsilon} f_r {\rm d}r, \ \ \ {\rm d}t-a.e.  \ t\in[0,T).
$$
\end{lem}
\begin{lem}\label{lem:2.3}
Let $0<T\leqslant+\infty$ and $(f_t)_{t\in[0,T]}$ be a real valued $(\F_t)$-progressively measurable process such that $||f_.||_{\hcal^\infty}<+\infty$. Then we have
$$ \  f_t=\lim\limits_{\varepsilon\To 0^+} \mathbf{E}\bigg[\frac{1}{\varepsilon}\int_t^{t+\varepsilon} f_r {\rm d}r\bigg| \F_t\bigg], \ \ \ {\rm d}t-a.e. \ t\in[0,T). $$
\end{lem}

The following \cref{lem:2.4} is taken from Lemma 3.1 in \cite{ZhangFan2013SPL}. \vspace{0.1cm}
\begin{lem}\label{lem:2.4}
Let random variables $(X_n)_{n\geqslant 1}$ and $X$ be square-integrable. If $\lim\limits_{n\to+\infty} X_n =X$ and $\sup\limits_{n\geqslant1}\mathbf{E}[|X_n|^2]<+\infty$, then  $$L^1-\lim\limits_{n\to +\infty}X_n=X.$$
\end{lem}

The following \cref{lem:2.5} gives a sequence of natural upper bounds for a quadratic function, which will play an important role in the proof of our main result. \vspace{0.1cm}
\begin{lem}\label{lem:2.5}
Let $f(y,z):\R\tim\R^d\To\R$. Assume that there exist constants $A,B\in\R_+$ such that for each $(y,z)\in\R\tim\R^d$, $|f(y,z)|\leqslant A+B|z|^2$. Then for each $n\in\mathbb{N}$, we have
\begin{equation}\label{eq:2.1}
\forall (y,z)\in\R\tim\R^d, \  |f(y,z)|\leqslant(n+B)(|y|+|z|^2)+\sup\limits_{\{(\bar{y},\bar{z}):|\bar{y}|+|\bar{z}|^2\leqslant \frac{A}{n}\}}|f(\bar{y},\bar{z})|.
\end{equation}
\end{lem}

\begin{proof}
If $|y|+|z|^2 \leqslant\frac{A}{n}$, then \eqref{eq:2.1} is clearly true. Otherwise, we have
$$|f(y,z)|\leqslant A+B|z|^2< n(|y|+|z|^2)+B|z|^2\leqslant(n+B)(|y|+|z|^2).$$
The proof of is complete.
\end{proof}

\section{Main result}
\label{sec:3-Main result}
\setcounter{equation}{0}

In this section, we will put forward and prove the general invariant representation for generators of general time interval quadratic BSDEs under assumptions \ref{A:H1}-\ref{A:H3} and some additional mild assumptions.

The following \cref{thm:3.1} is the main result of this paper.
\vspace{0.1cm}

\begin{thm}[Invariant Representation]\label{thm:3.1}
Assume that $0<T\leqslant+\infty$ and the generator $g$ satisfies assumptions {\rm \ref{A:H1}}-{\rm \ref{A:H3}}.

If the function $\varphi_\cdot$ in assumption {\rm \ref{A:H2}} also satisfies that for each $x\in\R_+$,
$$\int_0^T|\varphi_r(x)|^2{\rm d}r<+\infty\ \ {\rm and}\ \ \mathbf{E}[|\varphi_t(x)|^2]<+\infty, \ {\rm d}t-a.e. \ t\in[0,T),$$
then for each $(y,z)\in \R\tim\R^d$ and ${\rm d}t-a.e. \ t\in[0,T)$, we have
\begin{equation}\label{eq:3.1}
g(t,y,z)=L^1-\lim\limits_{\varepsilon\to 0^+}\frac{1}{\varepsilon}(Y_t^\varepsilon-y),
\end{equation}
where for each $\varepsilon\in (0,(T-t)\wedge 1]$, $(Y_t^\varepsilon,Z_t^\varepsilon)_{t\in[t,t+\varepsilon]}$ is any solution of the following BSDE (recalling \cref{pro:2.1}):
\begin{equation}\label{eq:3.2}
Y_s^{\varepsilon}=y+\langle z,B_{(t+\varepsilon)\wedge\tau}-B_t \rangle +\int_s^{t+\varepsilon}\mathbbm{1}_{r<\tau}g(r,Y_r^{\varepsilon},Z_r^{\varepsilon}){\rm d}r-\int_s^{t+\varepsilon}\langle Z_r^{\varepsilon},{\rm d}B_r\rangle, \ s\in[t,t+\varepsilon]
\end{equation}
with
\begin{equation*}
\tau:= {\rm inf}\bigg\{ s\geqslant t: |B_s-B_t|+\int_t^s|\varphi_r(K))|^2 {\rm d}r \geqslant 1\bigg\}
\end{equation*}
and
\begin{equation}\label{eq:3.3}
K:=3\bigg(|y|+|z|+\bigg\|\int_0^T f_r {\rm d}r\bigg\|_\infty \bigg) {\rm exp}\bigg( \bigg\|\int_0^T u_r {\rm d}r\bigg\|_\infty \bigg).\vspace{0.2cm}
\end{equation}

Furthermore, if the function $\varphi_\cdot$ in {\rm \ref{A:H2}} also satisfies  $\varphi_{\cdot}(x)\in\hcal^\infty(\R_+)$  for each $x\in\R_+$, then for each $(y,z)\in \R\times\R^d$ and ${\rm d}t-a.e. \ t \in [0,T)$, we have
\begin{equation}\label{eq:3.4}
 g(t,y,z)= \lim\limits_{\varepsilon\to 0^+}\frac{1}{\varepsilon}(Y_t^\varepsilon-y).
\end{equation}
\end{thm}

\begin{proof}
Fix $(t,y,z)\in [0,T)\tim\R\tim\R^d$ and choose $\varepsilon$ such that $0< \varepsilon \leqslant (T-t)\wedge1$. Assume that the generator $g$ satisfies assumptions {\rm \ref{A:H1}}-{\rm \ref{A:H3}}, and the function $\varphi_\cdot$ in \ref{A:H2} also satisfies that for each $x\in\R_+$,
$$\int_0^T|\varphi_r(x)|^2{\rm d}r<+\infty\ \ {\rm and}\ \ \mathbf{E}[|\varphi_t(x)|^2]<+\infty, \ {\rm d}t-a.e.\ t\in[0,T).$$
Suppose that $(Y_s^{\varepsilon},Z_s^{\varepsilon})_{s\in[t,t+\varepsilon]}$ is any solution of BSDE \eqref{eq:3.2} such that $Y_{\cdot}^\varepsilon$ is bounded and $\mathbf{E}[\int_t^{t+\varepsilon}|Z_s^\varepsilon|^2 {\rm d}s]<+\infty$. Then, $Y_s^\varepsilon\mathbbm{1}_{s\geqslant \tau}=y+\langle z, B_{(t+\varepsilon)\wedge\tau}-B_t\rangle$ and $Z_s^\varepsilon\mathbbm{1}_{s\geqslant \tau}=0$ for $s\in[t,t+\varepsilon]$. And, it follows from \cref{pro:2.1} and the definition of $\tau$ that
\begin{equation}\label{eq:3.5}
|{Y}_s^{\varepsilon}| \leqslant \bigg(|y|+|z|+\bigg\|\int_0^T f_r {\rm d}r\bigg\|_{\infty}\bigg) {\rm exp} \bigg(\bigg\| \int_0^T u_r {\rm d}r\bigg\|_{\infty}\bigg), \ \ \ s\in[t,t+\varepsilon].
\end{equation}

Now, we set, for each $s\in[t,t+\varepsilon]$,
$$\wid{Y}_s^{\varepsilon}:= Y_s^{\varepsilon}-y-\langle z, B_{s\wedge\tau}-B_t\rangle \ \ {\rm and} \ \ \wid{Z}_s^{\varepsilon}:= Z_s^{\varepsilon}-\mathbbm{1}_{s<\tau}z.$$
Then, it follows from \eqref{eq:3.5} that for each $s\in[t,t+\varepsilon]$,
\begin{equation}\label{eq:3.6}
|\wid{Y}_s^{\varepsilon}|\leqslant |Y_s^{\varepsilon}|+|y|+|z|\leqslant 2\bigg( |y|+|z|+\bigg\|\int_0^T f_r {\rm d}r\bigg\|_\infty \bigg){\rm exp}\bigg(\bigg\|\int_0^T u_r {\rm d}r \bigg\|_\infty\bigg)=: k.
\end{equation}
Let $q_k(\wid{y}) := \frac{k\wid{y}}{|\wid{y}|\vee k}$ for each $\wid{y}\in \R$. Then $|q_k(\wid{y})|\leqslant k \wedge|\wid{y}|$. By virtue of $\rm It\hat{o}$'s formula and the fact of $q_k(\wid{Y}_s^\varepsilon)= \wid{Y}_s^\varepsilon$ due to \eqref{eq:3.6}, it is not hard to verify that $(\wid{Y}_s^{\varepsilon}, \wid{Z}_s^{\varepsilon})_{s\in[t,t+\varepsilon]}$ is a solution of the following BSDE
\begin{equation}\label{eq:3.7}
\wid{Y}_s^{\varepsilon}=\int_s^{t+\varepsilon}\wid{g}(r,\wid{Y}_r^{\varepsilon},\wid{Z}_r^{\varepsilon}){\rm d}r-\int_s^{t+\varepsilon}\langle \wid{Z}_r^{\varepsilon},{\rm d}B_r \rangle,\ s\in[t,t+\varepsilon],
\end{equation}
where for each $(\wid{y}, \wid{z})\in\R\tim\R^d$,
$$\wid{g}(s,\wid{y},\wid{z}) := \mathbbm{1}_{s<\tau}g(s,q_k(\wid{y})+y+\langle z, B_{s \wedge \tau}-B_t\rangle, \wid{z}+z), \ s\in[t,t+\varepsilon].$$
By \ref{A:H2} together with the definition of $\tau$, ${\rm d}\mathbb{P}\times{\rm d} s-a.e.$, for each $(\wid{y},\wid{z})\in\R\times \R^d,$ we have
\begin{equation}\label{eq:3.8}
\begin{aligned}
&|\wid{g}(s,\wid{y},\wid{z})|\\
=&|\mathbbm{1}_{s<\tau} g(s,q_k(\wid{y})+y+\langle z, B_{s \wedge \tau}-B_t\rangle, \wid{z}+z)|\\
\leqslant &\mathbbm{1}_{s<\tau}\varphi_s(|q_k(\wid{y})+y+\langle z, B_{s \wedge \tau_\varepsilon}-B_t\rangle|) +h(|q_k(\wid{y})+y+\langle z, B_{s \wedge \tau}-B_t\rangle|)|\wid{z}+z|^2 \\
\leqslant & \mathbbm{1}_{s<\tau}\varphi_s(k+|y|+|z|)+2h(k+|y|+|z|)|\wid{z}|^2+2h(k+|y|+|z|)|z|^2\\
\leqslant &A(s)+B|\wid{z}|^2,
\end{aligned}
\end{equation}
 where $$A(s):= \mathbbm{1}_{s<\tau}\varphi_s(K)+2h(K)|z|^2 \ \ {\rm and} \ \ B:=2h(K)$$ with $K$ being defined in \eqref{eq:3.3}. It is clear that the process $A_\cdot$ satisfies that $\int_t^{t+\varepsilon}|A(r)|^2{\rm d}r<+\infty$ due to the assumptions of $\varphi_.$. It follows from \eqref{eq:3.8} that $\wid{g}$ satisfies assumptions \ref{A:H1}-\ref{A:H3} with $f_s=A(s)$, $u_s\equiv 0$, $\gamma \equiv B$, $\varphi_s(\cdot)\equiv A(s)$ and $h(\cdot)\equiv B$ on the time interval $[t,t+\varepsilon]$, and ${\rm d}\mathbb{P}\times{\rm d} s-a.e.$, for each $(\wid{y},\wid{z})\in\R\times \R^d$, we have
\begin{equation}\label{eq:3.9}
 |\wid{g}(s,\wid{y},\wid{z})-\wid{g}(s,0,0)| \leqslant 2A(s)+B|\wid{z}|^2.
\end{equation}
Then, by \cref{lem:2.5} we know that ${\rm d\mathbbm{P}\times ds-a.e.}$, for each $ (\wid{y},\wid{z})\in\R\times \R^d$ and $n\in\N$,
\begin{equation}\label{eq:3.10}
 \ |\wid{g}(s,\wid{y},\wid{z})-\wid{g}(s,0,0)|\leqslant (n+B)(|\wid{y}|+|\wid{z}|^2)+ \sup\limits_{(\bar{y},\bar{z})\in\mathcal{G}_n(s)} |\wid{g}(s,\bar{y},\bar{z})-\wid{g}(s,0,0)|,
\end{equation}
where and hereafter, $$\mathcal{G}_n(s) := \bigg\{(\bar{y},\bar{z}):|\bar{y}|+|\bar{z}|^2\leqslant \frac{2A(s)}{n}\bigg\}, \ \ s\in[t,t+\varepsilon].$$

With respect to $\widetilde{Y}_\cdot^\varepsilon$ and $\widetilde{Z}_\cdot^\varepsilon$, we have the following assertion, whose proof is partially inspired by Lemma 3.1 in Zheng and Li \cite{ZhengLi2015SPL}.\vspace{0.2cm}

\begin{pro}\label{pro:3.2}
For each $t\in[0,T)$, we have
\begin{equation}\label{eq:3.11}
\lim\limits_{\varepsilon\To 0^+}\frac{1}{\varepsilon} \mathbf{E}\bigg[\int_t^{t+\varepsilon}|\wid{Y}^{\varepsilon}_r| {\rm d}r \bigg|\F_t\bigg]=0,
\end{equation}

\begin{equation}\label{eq:3.12}
\lim\limits_{\varepsilon\To 0^+}\frac{1}{\varepsilon} \mathbf{E}\bigg[\int_t^{t+\varepsilon}|\wid{Z}^{\varepsilon}_r|^2 {\rm d}r \bigg]=0
\end{equation}
and
\begin{equation}\label{eq:3.13}
\lim\limits_{\varepsilon\To 0^+}\frac{1}{\varepsilon} \mathbf{E}\bigg[\int_t^{t+\varepsilon}|\wid{Z}^{\varepsilon}_r|^2 {\rm d}r \bigg|\F_t\bigg]=0.
\end{equation}
\end{pro}
\vspace{0.2cm}

\begin{proof}[Proof of \cref{pro:3.2}.]
Let $t\in[0,T)$ and  $\varepsilon\in(0,(T-t)\wedge1]$. Note that $\wid{g}$ satisfies \ref{A:H1}-\ref{A:H3} with $f_s=A(s)$, $u_s=0$ and $\gamma=B$.
It follows from \cref{pro:2.1} that BSDE \eqref{eq:3.7} admits a bounded solution $(\widetilde{Y}_s^\varepsilon,\widetilde{Z}_s^\varepsilon)_{s\in[t,t+\varepsilon]}$ and by H${\rm \ddot{o}}$lder's inequality together with the definitions of $A(\cdot)$ and $\tau$ as well as the fact of $0<\varepsilon<1$,
\begin{equation}\label{eq:3.14}
\begin{aligned}
 \sup\limits_{t\leqslant s \leqslant t+\varepsilon} |\widetilde{Y}_s^{\varepsilon}|
&\leqslant \bigg\|\int_t^{t+\varepsilon} A(s) {\rm d}s\bigg\|_{\infty} \leqslant \bigg\|\sqrt{\varepsilon}\bigg(\int_t^{t+\varepsilon} |A(s)|^2 {\rm d}s\bigg)^{\frac{1}{2}}\bigg\|_{\infty}
\leqslant \sqrt{\varepsilon} C
\end{aligned}
\end{equation}
with $$C:=\sqrt{2+8h^2(K)|z|^2},$$ which yields \eqref{eq:3.11} immediately.
Now we select $\varepsilon$ small enough such that
\begin{equation}\label{eq:3.15}
\sup\limits_{t \leqslant s \leqslant t+\varepsilon}|\widetilde{Y}_s^{\varepsilon}| \leqslant \frac{1}{4B}.
\end{equation}
Applying ${\rm It\hat{o}}$'s formula to $|\widetilde{Y}_s^{\varepsilon}|^2$ for $s\in[t,t+\varepsilon]$ and using \eqref{eq:3.8}, we have
\begin{equation}\label{eq:3.16}
\begin{aligned}
&|\widetilde{Y}_t^{\varepsilon}|^2 +\int_t^{t+\varepsilon}|\widetilde{Z}_r^{\varepsilon}|^2 {\rm d}r\\
\leqslant &2\int_t^{t+\varepsilon} |\widetilde{Y}_r^{\varepsilon}|{\rm sgn} (\widetilde{Y}_r^{\varepsilon}) \widetilde{g}(r,\widetilde{Y}_r^{\varepsilon},\widetilde{Z}_r^{\varepsilon}) {\rm d}r-2\int_t^{t+\varepsilon} \langle \widetilde{Y}_r^{\varepsilon}, \widetilde{Z}_r^{\varepsilon} {\rm d}B_r\rangle\\
\leqslant &2\int_t^{t+\varepsilon} A(r)|Y_r^{\varepsilon}|{\rm d}r + 2\int_t^{t+\varepsilon} B|\wid{Y}_r^{\varepsilon}||\widetilde{Z}_r^{\varepsilon}|^2{\rm d}r-2\int_t^{t+\varepsilon} \langle \wid{Y}_r^{\varepsilon}, \wid{Z}_r^{\varepsilon}{\rm d}B_r\rangle.
\end{aligned}
\end{equation}
Note that $\big(\int_0^t \langle \widetilde{Y}_r^\varepsilon,\widetilde{Z}_r^\varepsilon {\rm d}B_r\rangle\big)_{t\in[0,T]}$ is a martingale starting from 0. It follows from \eqref{eq:3.14}, \eqref{eq:3.15}, \eqref{eq:3.16} and ${\rm H\ddot{o}lder}$'s inequality that for each $0\leqslant u \leqslant t$,
\begin{equation*}
\begin{aligned}
\frac{1}{\varepsilon} \mathbf{E}\bigg[\int_t^{t+\varepsilon}|\wid{Z}_r^{\varepsilon}|^2 {\rm d}r\bigg|\F_u\bigg]
&\leqslant \frac{2}{\varepsilon}\mathbf{E} \bigg[\int_t^{t+\varepsilon} A(r)|\wid{Y}_r^{\varepsilon}| {\rm d}r\bigg|\F_u\bigg]\\
&\leqslant \frac{2}{\varepsilon} \mathbf{E} \bigg[ \sqrt{\varepsilon} C \bigg\| \sqrt{\varepsilon} \bigg(\int_t^{t+\varepsilon}|A(r)|^2 dr\bigg)^{\frac{1}{2}} \bigg\|_{\infty}\bigg|\F_u\bigg] \\
&\leqslant 2C\bigg\|\int_t^{t+\varepsilon} |A(r)|^2 {\rm d}r\bigg\|_\infty^{\frac{1}{2}},
\end{aligned}
\end{equation*}
which yields \eqref{eq:3.12} and \eqref{eq:3.13}. \cref{pro:3.2} is then proved.
\end{proof}
\vspace{0.1cm}

In the sequel, letting $s=t$ and taking the conditional expectation with respect to $\F_t$ in both sides of BSDE \eqref{eq:3.7} leads to the following identity
\begin{equation*}
 \frac{1}{\varepsilon}(Y_t^{\varepsilon}-y)=\frac{1}{\varepsilon}\wid{Y}_t^{\varepsilon} =\frac{1}{\varepsilon} \mathbf{E} \bigg[\int_t^{t+\varepsilon}\wid{g}(r,\wid{Y}_r^{\varepsilon},\wid{Z}_r^{\varepsilon}){\rm d}r\bigg|\F_t\bigg].
\end{equation*}
Set
$$
M_t^{\varepsilon}:= \frac{1}{\varepsilon} \mathbf{E} \bigg[\int_t^{t+\varepsilon}\wid{g}(r,\wid{Y}_r^{\varepsilon},\wid{Z}_r^{\varepsilon}){\rm d}r\bigg|\F_t\bigg] \ \ \  {\rm and} \ \ \
N_t^{\varepsilon}:= \frac{1}{\varepsilon} \mathbf{E}\bigg[\int_t^{t+\varepsilon}\wid{g}(r,0,0){\rm d}r\bigg|\F_t\bigg].
$$
It follows that for each $t\in[0,T)$,
\begin{equation}\label{eq:3.17}
\frac{1}{\varepsilon}(Y_t^{\varepsilon}-y)-g(t,y,z)= \frac{1}{\varepsilon}\wid{Y}_t^{\varepsilon}- \wid{g}(t,0,0)=(M_t^{\varepsilon}-N_t^{\varepsilon})+(N_t^{\varepsilon}-\wid{g}(t,0,0)).
\end{equation}
Then, in order to prove \eqref{eq:3.1}, it is sufficient to prove that for ${\rm d}t-a.e.\ t\in[0,T)$, $(M_t^{\varepsilon}-N_t^{\varepsilon})$ and $(N_t^{\varepsilon}-\wid{g}(t,0,0))$ respectively tend to 0 in $L^1$ sense as $\varepsilon\To 0^+$.\vspace{0.2cm}

We first treat the term $(M_t^{\varepsilon}-N_t^{\varepsilon})$. By Jensen's inequality and \eqref{eq:3.10}, we get that for each $ t\in[0,T)$ and $n\in\N$,
\begin{equation}\label{eq:3.18}
\begin{aligned}
\mathbf{E}[|M_t^{\varepsilon}-N_t^{\varepsilon}|]&=\mathbf{E}\bigg\{\bigg| \frac{1}{\varepsilon}\mathbf{E} \bigg[\int_t^{t+\varepsilon} (\wid{g}(r,\wid{Y}_r^{\varepsilon},\wid{Z}_r^{\varepsilon}) -\wid{g}(r,0,0)){\rm d}r\bigg|\F_t\bigg]\bigg|\bigg\}\\
&\leqslant \mathbf{E}\bigg[\frac{1}{\varepsilon} \int_t^{t+\varepsilon}|\wid{g}(r,\wid{Y}_r^{\varepsilon},\wid{Z}_r^{\varepsilon}) -\wid{g}(r,0,0)|{\rm d}r \bigg] \\
&\leqslant   (n+B)\mathbf{E}\bigg[\frac{1}{\varepsilon}\int_t^{t+\varepsilon}| \wid{Y}_r^{\varepsilon}|{\rm d}r \bigg] +(n+B)\mathbf{E}\bigg[\frac{1}{\varepsilon} \int_t^{t+\varepsilon}|\wid{Z}_r^{\varepsilon}|^2{\rm d}r \bigg]\\
&\ \ \ \ +\mathbf{E}\bigg[  \frac{1}{\varepsilon} \int_t^{t+\varepsilon} \sup\limits_{(\bar{y},\bar{z})\in\mathcal{G}_n(r)} |\wid{g}(r,\bar{y},\bar{z})-\wid{g}(r,0,0)|{\rm d}r \bigg].
\end{aligned}
\end{equation}\vspace{0.1cm}
It follows from \eqref{eq:3.9} and the definition of $\mathcal{G}_n(r)$ that for each $n\in\N$ and ${\rm d}r-a.e. \ r\in[t,t+\varepsilon]$, we have
\begin{equation}\label{eq:3.19}
\begin{aligned}
\sup\limits_{(\bar{y},\bar{z})\in\mathcal{G}_n(r)} |\wid{g}(r,\bar{y},\bar{z})-\wid{g}(r,0,0)|  &\leqslant \sup\limits_{(\bar{y},\bar{z})\in\mathcal{G}_n(r)} \big(2A(r)+B|\bar{z}|^2\big) \\
&\leqslant 2A(r)+B\frac{2A(r)}{n}\leqslant 2(B+1)A(r).
\end{aligned}
\end{equation}
Note that $\int_0^{T\wedge m}A(t){\rm d}t<+\infty$ for each $m\in\N$. Thanks to \cref{lem:2.2}, we have that ${\rm d}t-a.e. \ t\in [0,T)$,
\begin{equation}\label{eq:3.20}
\begin{aligned}
\lim\limits_{\varepsilon\to0^+} \frac{1}{\varepsilon}\int_t^{t+\varepsilon} \sup\limits_{(\bar{y},\bar{z})\in\mathcal{G}_n(r)} |\wid{g}(r,\bar{y},\bar{z})-\wid{g}(r,0,0)|{\rm d}r= \sup\limits_{(\bar{y},\bar{z})\in\mathcal{G}_n(t)} |\wid{g}(t,\bar{y},\bar{z})-\wid{g}(t,0,0)|.
\end{aligned}
\end{equation}
And, in view of \eqref{eq:3.19}, applying H${\rm \ddot{o}}$lder's inequality yields that ${\rm d}t-a.e. \ t\in[0,T)$,
\begin{equation}\label{eq:3.21}
\begin{aligned}
&\sup\limits_{\varepsilon>0}\mathbf{E}\bigg[\bigg(\frac{1}{\varepsilon} \int_t^{t+\varepsilon}\sup\limits_{(\bar{y},\bar{z})\in\mathcal{G}_n(r)} |\wid{g}(r,\bar{y},\bar{z})-\wid{g}(r,0,0)|{\rm d}r\bigg)^2\bigg]  \\
\leqslant &4(B+1)^2 \sup\limits_{\varepsilon>0} \mathbf{E} \bigg[\frac{1}{\varepsilon^2}\bigg(\int_t^{t+\varepsilon}A(r){\rm d}r\bigg)^2\bigg]\\
\leqslant &(B+1)^2 \sup\limits_{\varepsilon>0} \mathbf{E}\bigg[\frac{1}{\varepsilon}\int_t^{t+\varepsilon} |A(r)|^2 {\rm d}r\bigg] <+\infty.
\end{aligned}
\end{equation}
In the last inequality we have used the following fact, due to Fubini's theorem and \cref{lem:2.2},
\begin{equation}\label{eq:3.22}
\lim\limits_{\varepsilon\To 0^+} \mathbf{E}\bigg[\frac{1}{\varepsilon}\int_t^{t+\varepsilon} |A(r)|^2 {\rm d}r \bigg] =\mathbf{E}\bigg[ |A(t)|^2 \bigg], \ {\rm d}t-a.e. \ t\in[0,T).
\end{equation}
Then, in view of \eqref{eq:3.20} and \eqref{eq:3.21}, it follows from \cref{lem:2.4} that for each $n\in\N$ and ${\rm d}t-a.e. \ t\in[0,T)$, we have
\begin{equation}\label{eq:3.23}
\begin{aligned}
&\lim\limits_{\varepsilon\To 0^+}  \mathbf{E}\bigg[ \frac{1}{\varepsilon} \int_t^{t+\varepsilon} \sup\limits_{(\bar{y},\bar{z})\in\mathcal{G}_n(r)} |\wid{g}(r,\bar{y},\bar{z})-\wid{g}(r,0,0)| {\rm d}r \bigg]\\
=&\mathbf{E}\bigg[\sup\limits_{(\bar{y},\bar{z})\in\mathcal{G}_n(t)} |\wid{g}(t,\bar{y},\bar{z}) -\wid{g}(t,0,0)|\bigg].
\end{aligned}
\end{equation}
Furthermore, it follows from \eqref{eq:3.19}, Lebesgue's dominated convergence, the definition of $\mathcal{G}_n(t)$ and the continuity of $g$ in $(y,z)$ that for ${\rm d}t-a.e. \ t\in[0,T),$
\begin{equation}\label{eq:3.24}
\begin{aligned}
&\lim\limits_{n\To \infty}  \mathbf{E}\bigg[ \sup\limits_{(\bar{y},\bar{z})\in\mathcal{G}_n(t)} |\wid{g}(t,\bar{y},\bar{z})-\wid{g}(t,0,0)|\bigg]\\
=&\mathbf{E}\bigg[\lim\limits_{n\To \infty} \sup\limits_{(\bar{y},\bar{z})\in\mathcal{G}_n(t)} |\wid{g}(t,\bar{y},\bar{z})-\wid{g}(t,0,0)|\bigg]=0.
\end{aligned}
\end{equation}
Thus, in view of \eqref{eq:3.11}, \eqref{eq:3.12}, \eqref{eq:3.23} and \eqref{eq:3.24}, sending first $\varepsilon\To 0^+$ and then $n\To \infty$ in \eqref{eq:3.18} yields that for ${\rm d}t-a.e.\ t\in[0,T)$,
\begin{equation}\label{eq:3.25}
\lim\limits_{\varepsilon\To 0^+}\mathbf{E}[|M_t^{\varepsilon}-N_t^{\varepsilon}|]
=0.
\end{equation}

We now consider the term $(N_t^{\varepsilon}-\wid{g}(t,0,0))$. By \eqref{eq:3.8}, H${\rm \ddot{o}}$lder's inequality and \eqref{eq:3.22}, we obtain that
\begin{equation*}
\begin{aligned}
&\sup\limits_{\varepsilon>0} \mathbf{E}\bigg[ \bigg(\frac{1}{\varepsilon} \int_t^{t+\varepsilon} |\wid{g}(r,0,0)-\widetilde{g}(t,0,0)| {\rm d}r\bigg)^2 \bigg]\\ \leqslant &2\sup\limits_{\varepsilon>0} \mathbf{E} \bigg[\frac{1}{\varepsilon} \int_t^{t+\varepsilon} |A(r)|^2 {\rm d}r\bigg]+2\mathbf{E}\big[|A(t)|^2]<+\infty, \ \ {\rm d}t-a.e. \ t\in[0,T).
\end{aligned}
\end{equation*}
Then, it follows from \cref{lem:2.2} and \cref{lem:2.4} that for ${\rm d}t-a.e.\ t\in[0,T)$,
\begin{equation}\label{eq:3.26}
\begin{aligned}
\lim\limits_{\varepsilon\To 0^+}\mathbf{E}[|N_t^{\varepsilon}-\wid{g}(t,0,0)|] &\leqslant \lim\limits_{\varepsilon\To 0^+} \mathbf{E}\bigg[\frac{1}{\varepsilon}\int_t^{t+\varepsilon} |\wid{g}(r,0,0)-\wid{g}(t,0,0)| {\rm d}r\bigg]\\
&= \mathbf{E}\bigg[\lim\limits_{\varepsilon\To 0^+} \frac{1}{\varepsilon}\int_t^{t+\varepsilon}|\wid{g}(r,0,0)-\wid{g}(t,0,0)| {\rm d}r\bigg] = 0.
\end{aligned}
\end{equation}
Thus, the desired equality \eqref{eq:3.1} follows from \eqref{eq:3.17}, \eqref{eq:3.25} and \eqref{eq:3.26}.\vspace{0.2cm}

Finally, we further suppose that the function $\varphi_\cdot$ in \ref{A:H2} also satisfies  $\varphi_{\cdot}(x)\in\hcal^\infty(\R_+)$ for each $x\in\R_+$. Then, $A(\cdot)\in\hcal^\infty(\R_+)$, and from \eqref{eq:3.17} we know that in order to prove \eqref{eq:3.2}, it is sufficient to prove that for ${\rm d}t-a.e.\ t\in[0,T)$, $(M_t^{\varepsilon}-N_t^{\varepsilon})$ and $(N_t^{\varepsilon}-\wid{g}(t,0,0))$ respectively tend to 0 in $\ps$ sense as $\varepsilon\To0^+$.\vspace{0.2cm}

We first treat the term $(M_t^{\varepsilon}-N_t^{\varepsilon})$. It follows from \eqref{eq:3.10} that for each $t\in [0,T)$ and $n\in \N$, we have
\begin{equation}\label{eq:3.27}
\begin{aligned}
|M_t^{\varepsilon}-N_t^{\varepsilon}|&\leqslant \mathbf{E}\bigg[\frac{1}{\varepsilon} \int_t^{t+\varepsilon}|\wid{g}(r,\wid{Y}_r^{\varepsilon},\wid{Z}_r^{\varepsilon}) -\wid{g}(r,0,0)|{\rm d}r \bigg|\F_t\bigg] \\
&\leqslant   (n+B)\mathbf{E}\bigg[\frac{1}{\varepsilon}\int_t^{t+\varepsilon} |\wid{Y}_r^{\varepsilon}|{\rm d}r \bigg|\F_t\bigg] +(n+B)\mathbf{E}\bigg[\frac{1}{\varepsilon}\int_t^{t+\varepsilon} |\wid{Z}_r^{\varepsilon}|^2{\rm d}r \bigg|\F_t\bigg]\\
&\ \ +\mathbf{E}\bigg[ \frac{1}{\varepsilon} \int_t^{t+\varepsilon} \sup\limits_{(\bar{y},\bar{z})\in\mathcal{G}_n(r)} |\wid{g}(r,\bar{y},\bar{z})-\wid{g}(r,0,0)|{\rm d}r \bigg|\F_t\bigg].
\end{aligned}
\end{equation}
Furthermore, by \eqref{eq:3.9} and the definition of $\mathcal{G}_n(r)$ we know that for each $n\in\N$ and ${\rm d}r-a.e. \ r\in[t,t+\varepsilon]$,
$$\sup\limits_{(\bar{y},\bar{z})\in\mathcal{G}_n(r)} |\wid{g}(r,\bar{y},\bar{z})-\wid{g}(r,0,0)| \leqslant \sup\limits_{(\bar{y},\bar{z})\in\mathcal{G}_n(r)} (2A(r)+B|\bar{z}|^2) \leqslant 2(B+1)\|A(\cdot)\|_{\hcal^\infty}.$$
Then, applying \cref{lem:2.3} yields that for each $n\in\N$ and ${\rm d}t-a.e.\ t\in[0,T)$,
\begin{equation}\label{eq:3.28}
\begin{aligned}
&\lim\limits_{\varepsilon\To 0^+}\mathbf{E}\bigg[\frac{1}{\varepsilon} \int_t^{t+\varepsilon}\sup\limits_{(\bar{y},\bar{z})\in\mathcal{G}_n(r)} |\wid{g}(r,\bar{y},\bar{z})-\wid{g}(r,0,0)|{\rm d}r\bigg|\F_t\bigg]\\
=&\sup\limits_{(\bar{y},\bar{z})\in\mathcal{G}_n(t)} |\wid{g}(t,\bar{y},\bar{z})-\wid{g}(t,0,0)|.
\end{aligned}
\end{equation}
Thus, in view of \eqref{eq:3.11}, \eqref{eq:3.13}, \eqref{eq:3.28}, the definition of $\mathcal{G}_n(t)$ and the continuity of $g$ in $(y,z)$, sending first $\varepsilon\To 0^+$ and then $n\To \infty$ in \eqref{eq:3.27} yields that for ${\rm d}t-a.e. \ t\in[0,T)$
\begin{equation}\label{eq:3.29}
\lim\limits_{\varepsilon\To 0^+}|M_t^{\varepsilon}-N_t^{\varepsilon}|=0.
\end{equation}

Now we consider the term $(N_t^{\varepsilon}-\wid{g}(t,0,0))$. Note from \eqref{eq:3.8} that ${\rm d}t-a.e.\ t\in[0,T)$,
$$|\wid{g}(t,0,0)|\leqslant A(t) \leqslant \|A(\cdot)\|_{\hcal^\infty}.$$
Applying \cref{lem:2.3} yields that for ${\rm d}t-a.e.\ t\in[0,T)$,
\begin{equation}\label{eq:3.30}
\lim\limits_{\varepsilon\To 0^+}|N_t^{\varepsilon}-\wid{g}(t,0,0)|\leqslant \lim\limits_{\varepsilon\To 0^+}\mathbf{E}\bigg[\frac{1}{\varepsilon}\int_t^{t+\varepsilon}|\wid{g}(r,0,0)-\wid{g}(t,0,0)| {\rm d}r\bigg|\F_t\bigg] = 0
\end{equation}
Thus, the desired assertion \eqref{eq:3.4} follows from \eqref{eq:3.17}, \eqref{eq:3.29} and \eqref{eq:3.30}.

The proof of \cref{thm:3.1} is then complete.\vspace{0.1cm}
\end{proof}

We end this paper by the following Remark.\vspace{0.1cm}

\begin{rmk}\label{rmk:3.3}
(i) It following from \cref{pro:2.1} that assumptions \ref{A:H1}-\ref{A:H3} of the generator $g$ guarantee the existence of bounded solutions of BSDEs. However, in order to obtain the representation of the generator, some stronger integrability conditions on the function $\varphi_\cdot$ in \ref{A:H2} are additionally required, which is the same as usual. It is not vary hard to find the reason from the proof of \cref{thm:3.1}.

(ii) We would like to mention that the conditions of the generator $g$ in \cref{thm:3.1} are strictly weaker than all the corresponding representation results mentioned in the introduction. In fact, \cref{thm:3.1} strengthens and unifies these known results.

(iii) In the proof of \cref{thm:3.1}, we systematically apply the stopping time technique, the truncation technique and the a priori estimate technique respectively used in \cite{MaYao2010SAA}, \cite{XiaoFan2017SPL}, \cite{ZhengLi2018AMASES} and \cite{LuoFan2017SD}. In particular, we also use \cref{lem:2.5} to obtain \eqref{eq:3.10} in the proof of \cref{thm:3.1}, which is a key point of the proof and very different from the proofs of existing results. The use of \cref{lem:2.5} replaces the convolution approximation technique used in \cite{FanJiang2010JCAM}, \cite{FanJiang2013AMS}, \cite{XiaoFan2017SPL}-\cite{ZhengLi2018AMASES} and simplifies significantly the whole proof \cref{thm:3.1}, which is a natural and intrinsic idea.
\end{rmk}

\vspace{0.1cm}



\end{document}